\documentclass[10pt]{article}
\usepackage{amsfonts,amssymb,amscd,amsmath,amsthm,enumerate,verbatim,calc,graphicx,latexsym,dsfont}
\usepackage[]{color}
\usepackage {oldgerm}
\usepackage{enumitem}
\setenumerate{label={\normalfont(\arabic*)}, ref={\normalfont(\arabic*)}, itemsep=0pt}

\newtheorem{theorem}{Theorem}
\newtheorem{lemma}{Lemma}

\newtheorem{question}[theorem]{Question}

\theoremstyle{definition}

\newcommand{\A}{\ensuremath{\operatorname{Aut}}}
\newcommand{\Sym}{\ensuremath{\operatorname{Sym}}}

\newcommand{\m}{\ensuremath{\operatorname{m}}}

\newcommand{\id}{\ensuremath{\text{\rm id}}}

\newcommand{\conf}{\ensuremath{\operatorname{conf}}}

\newcommand{\N}{{\mathbb{N}}}

\renewcommand{\proof}{\noindent{\bf Proof\,\ }}
\renewcommand{\qed}{\hfill $\Box$\medskip}

\begin{document}
\title
{On a theorem of Halin\\
\vspace{4mm}\Large Dedicated to the memory of Rudolf Halin}
\author{
Wilfried Imrich\\
Department Mathematics and Information Technology,\\
Montanuniversit\"{a}t Leoben, 8700 Leoben, Austria \\
imrich@unileoben.ac.at
\and
Simon M.~Smith\\
Department of Mathematics,
NYC College of Technology,\\
City University of New York,
New York, NY, USA\\
sismith@citytech.cuny.edu
}
\date{\today}
\maketitle

\begin{abstract}
This note presents a new, elementary proof of a generalization of a theorem of Halin to graphs with unbounded degrees, which is then applied to show that every connected, countably infinite graph $G$, with $\aleph_0 \leq |\A(G)| < 2^{\aleph_0}$ and subdegree-finite automorphism group, has a finite set $F$ of vertices that is setwise stabilized only by the identity automorphism. A bound on the size of such sets,  which are called {\em distinguishing},  is also provided.

To put this  theorem of Halin and its generalization into perspective, we also discuss several related non-elementary, independent results and their methods of proof.
\end{abstract}

\vspace{0.2cm}


\section{Introduction}
\label{sec:intro}
\setcounter{footnote}{1}
This  is an homage to Halin and  his  seminal paper  from 1973 on automorphisms and endomorphisms of infinite, locally finite graphs \cite{ha-73}.
It  focusses on Theorem 6 of that paper, which asserts that a locally finite, connected graph  with infinite automorphism group has a finite base if and only if its automorphism group  is countable.

Halin and most graphs theorists were seemingly not aware of the fact that  stronger forms of this theorem, which did not need the assumption of local finiteness, had already been published in 1968 in the language of infinitary languages by Kueker (\cite[Theorem 2.1]{kueker}) and 1970 in the language of mathematical logic by Reyes (\cite[Theorem 2.2.2]{Reyes}). In 1987, Evans (\cite[Theorem 1.1]{Evans87}) published a version more accessible to graph theorists in Archiv der Mathematik. It implies that Halin's Theorem 6 from \cite{ha-73}  holds without the assumption of local finiteness.  We formulate this as Theorem \ref{thm:halin1} and
derive it with an elementary argument inspired by Halin's original proof.

We then apply this result to show that every graph  with a finite base and infinite automorphism group that is subdegree
finite\footnote{A permutation group is subdegree-finite if its point stabilizers have finite orbits.} contains a set
that is (setwise) stabilized only by the identity automorphism\footnote{Such sets are called \emph{distinguishing sets}.}
and that the  size of this set is less than three times the size of the base, see Theorem \ref{thm:Cameron}. This generalizes \cite[Theorem 3.2]{boim-2015} from locally finite graphs to graphs of larger valence.

In the last section  we discuss Evans' and other approaches in more detail.
These approaches take recourse to results about topological and metric spaces, including a theorem of Baire, and usually yield more general results.

We conclude with  the problem of extending  Halin's theorem to uncountable graphs.

\section{Preliminaries}

We focus on connected, countably infinite graphs $G$   with unbounded degrees and are interested in sets of vertices that are only (setwise) stabilized by the trivial automorphism. Such sets are called \emph{distinguishing}. Not every graph has such a set, but if it does, then we call the minimum size of such a set the \emph{distinguishing cost $\rho(G)$} of $G$. There are large classes of infinite graphs that have such sets, and often $\rho(G)$ is finite. These are the graphs we are interested in.

The distinguishing cost was introduced and studied by Boutin \cite{bo-2008} for finite graphs. Infinite graphs with finite cost were treated for the first time in \cite{boim-2015}. These investigations did not cover graphs with infinite degrees, that is, they were restricted to \emph{locally finite} graphs.

Recall that the \emph{set stabilizer} of $S \subseteq V(G)$, denoted $\A(G)_S$, is the set of all  automorphisms $\varphi$ for which
\[\varphi(x) \in S \iff x\in S.\]
We also  say that $S$ is \emph{invariant} under $\varphi$, or that $\varphi$ \emph{preserves} $S$, and  write $\varphi(S)=S$.
The \emph{point stabilizer} of $S$, denoted by $\A(G)_{(S)}$
 is the set of all $\varphi\in \A(G)$ for which $\varphi(x)=x$ for all $x\in S$.
 If $A$ is a set of automorphisms and $x \in V(G)$, then the \emph{orbit} $A(x)$  of $x$ under the action of $A$ is the set $\{\varphi(x) : \varphi \in A\}$.

  A \emph{base} of $\A(G)$   is a set of vertices $S$ whose point stabilizer is trivial, and the \emph{ determining number} of the graph $G$, denoted  by ${\det}(G)$, is the minimum size of a base of $G$.
  Clearly every distinguishing set also is a base. Theorem \ref{thm:Cameron} below estimates the cost as a function of the size of  $\det(G)$. It generalizes \cite[Theorem  3.2]{boim-2015} for locally finite graphs.

  Notice that every base
  $S$ has the property that whenever $\varphi, \psi \in {\A}(G)$ so that $\varphi(x)=\psi(x)$ for all $x\in {S}$, then $\varphi=\psi$. Thus every automorphism of $G$ is uniquely determined by its action on the vertices of a base.

Furthermore, for $F\subseteq V(G)$, an automorphism $\varphi\in \A(G)_F$ can be thought of as a permutation in $\Sym(F)$ by restricting the action of $\varphi$ to $F$, denoted $\varphi|F$.  Thus we have a natural map $\A(G)_F \to \A(G)|F \leq \Sym(F)$.  Note that this map is injective if and only if $F$ is a base for $G$.  In this case $\A(G)_F \cong \A(G)|F$.   If $F$ is a finite base for $\A(G)$, then $\A(G)_F$ is a finite group.

The class of connected, locally finite, infinite graphs is usually denoted by $\Gamma$. Graphs in $\Gamma$ have the property that the orbits of all point stabilizers are finite. Automorphism groups with this property are called \emph{subdegree-finite}.
We  do not require subdegree-finiteness for Theorem \ref{thm:halin1}, but for Theorem \ref{thm:Cameron} it is essential.

The {\it motion} of an automorphism $\varphi\in \A(G)$, denoted $\m(\varphi)$,  is the number of vertices moved by $\varphi$.  The {\it motion} of the graph $G$, denoted by $\m(G)$, is the minimum motion of the nontrivial elements of $\A(G)$.

The vertices of a graph with a distinguishing set obviously admit a 2-coloring that is only preserved by the identity automorphism. Such a coloring is called a 
\emph{distinguishing 2-coloring}.
These colorings are a special case of distinguishing \emph{$d$-colorings}, where  the vertices of a graph $G$ are labeled with the integers $1,\ldots, d$ in such a way that only the trivial automorphism preserves the colors.  A graph with such a coloring is called \emph{$d$-distinguishable}.  This concept was introduced by Albertson and Collins in \cite{alco-96} and initiated a wealth of results for both finite and infinite graphs. For infinite graphs, \cite{istw-2015} gives a good introduction with numerous further references to finite and infinite graphs. But, we also wish to point out a paper of Babai \cite{ba-1977} from 1977, in which he proves,  in a different setting, deep results about distinguishing infinite graphs.

\section{Halin's theorem for graphs with unbounded degrees}

In this section we extend  Halin's theorem, which is listed below,  to countable, connected graphs that are not in $\Gamma$.

\begin{theorem}\label{thm:halin} {\rm \cite{ha-73} (Halin 1973)} A connected, locally finite, infinite graph $G$ has uncountable $\A(G)$ if and only if for every finite $F \subset V(G)$ there exists a nontrivial automorphism $\varphi$ of $G$ such that $\varphi(v) = v$ for each $v \in F$.\end{theorem}

In other words, this theorem says that a graph $G \in \Gamma$ has countable automorphism group if and only if it  has a finite base.  We are interested in
graphs with infinite automorphism group and prove the following extension to Halin's Theorem, which is the graph theoretic version of Kueker \cite[Theorem 2.1]{kueker} and Evans \cite[Theorem 1.1]{Evans87}.
Notice that Kueker's result precedes that of Halin, but not the one of Evans. Interestingly, Evans' paper also contains an outline of an elegant proof of Evans' Theorem suggested by Theo Grundh\"{o}fer which invokes Baire's Theorem (see Section~\ref{section:nonelementaryresults}).

We now give
the  extension of Theorem \ref{thm:halin} and
its
elementary 
combinatorial
proof.

\begin{theorem}\label{thm:halin1}Let $G$ be a graph with countably  many vertices. Then $\A(G)$ is countable or has cardinality $2^{\aleph_0}$, with $\A(G)$  countable if and only if $\A(G)$ has a finite base.
\end{theorem}

\proof The theorem is obviously true if $G$ is finite. Suppose $F \subseteq V(G)$. The set of images of $F$ under $\A(G)$ 
has cardinality
$|\A(G) : \A(G)_F| \leq \aleph_0$.
If $F$ is a finite base for $\A(G)$, then
$\A(G)_F$ is a finite group 
and $|\A(G)| = |\A(G) : \A(G)_F| \cdot |\A(G)_F : \A(G)_{(F)}|$. 
Hence $\A(G)$ is countable.

Assume that $\A(G)$  has no finite base, that is, to every finite set $F \subseteq V(G)$ there exists a nontrivial  automorphism of $G$ that fixes every element of $F$.
Fix some enumeration $\{v_0, v_1, v_2, \ldots\}$ of $V(G)$. Let $E$ be the set of all infinite sequences
$\{ \epsilon = (\epsilon_0, \epsilon_1, \epsilon_2, \ldots) : \epsilon_i \in \{0, 1\} \}$.

We will presently be defining inductively, for each integer $i \geq 0$, a finite set $F_i$, a vertex $x_i$, and an automorphism $\varphi_i$ which fixes each vertex in $F_i$ but does not fix $x_i$. To simplify our notation, for each
$\epsilon = (\epsilon_0, \epsilon_1, \ldots ) \in E$, let
$\alpha_i^\epsilon$ denote the automorphism $\varphi_0^{\epsilon_0} \cdots \varphi_i^{\epsilon_i}$
and
let $\alpha_i^{-\epsilon}$ denote the automorphism $\varphi_i^{-\epsilon_i} \cdots \varphi_0^{-\epsilon_0}$; write
$\alpha_i^E := \{\alpha_i^\epsilon : \epsilon \in E\}$ and
$\alpha_i^{-E} := \{\alpha_i^{-\epsilon} : \epsilon \in E\}$.
Note that $\alpha_i^E$ and $\alpha_i^{-E}$ are finite sets of automorphisms and each set contains the identity automorphism. We will write $\alpha_i^E(F_i)$ (resp. $\alpha_i^{-E}(F_i)$) to denote the set of all images of vertices in $F_i$ under automorphisms in $\alpha_i^E$ (resp. $\alpha_i^{-E}$). Note that $\alpha_i^{E}(F_i)$ and $\alpha_i^{-E}(F_i)$ are also finite sets.

Let us now begin our inductive construction. Take $F_0$ to be any finite set of vertices containing $v_0$ and let $\varphi_0$ be an automorphism that fixes every vertex of $F_0$, but which is  not the identity. Then there exists  a vertex $x_0$ that is moved by $\varphi_0$.

Let $F_1$ be a finite set of vertices containing $\alpha^E_0(F_0) \cup \alpha^{-E}_0(F_0) \cup \{x_0, v_1\}$. By assumption, there exists an automorphism $\varphi_1$ which fixes every vertex in $F_1$ but does not fix some vertex $x_1$.

Suppose, in addition, that for all integers $i$ satisfying $1 < i \leq k$ we have chosen a finite set $F_i$ which contains $\alpha_{i-1}^E(F_{i-1}) \cup \alpha_{i-1}^{-E}(F_{i-1}) \cup \{x_{i-1}, v_{i}\}$, and we have chosen a nontrivial automorphism $\varphi_i$ which fixes every vertex in $F_i$ but does not fix $x_i$. Let $F_{k+1}$ be a finite set of vertices containing
$\alpha_k^E(F_k) \cup \alpha_k^{-E}(F_k) \cup \{x_k, v_{k+1}\}$. By assumption, there exists an automorphism $\varphi_{k+1}$ which fixes every vertex in $F_{k+1}$ but does not fix some vertex $x_{k+1}$.

Thus, we have $v_k \in F_k$ for all $k \geq 0$, with $F_0 \subseteq F_1 \subseteq \cdots$ and $\bigcup_{k = 0}^\infty F_k = V(G)$.
If $\ell > k \geq 0$ and $\epsilon \in E$, then for all $v \in F_k$ we have $\alpha_k^{\pm \epsilon}(v) \in F_{k+1} \subseteq F_\ell$. Hence, one may easily verify that
\[\alpha_k^\epsilon(v) = \alpha_\ell^\epsilon(v)
\quad \text{and} \quad
\alpha_k^{-\epsilon}(v) = \alpha_\ell^{-\epsilon}(v).\]

We now define, for each $\epsilon \in E$, a map $\alpha^\epsilon : V(G) \rightarrow V(G)$ via
\[\alpha^\epsilon(v_k) := \alpha_k^\epsilon(v_k) \quad \quad \text{for all $k \geq 0$.}\]
Each map is well-defined on its domain $V(G) = \{v_0, v_1, \ldots\}$. For all $v_k, v_\ell \in V(G)$  choose an integer $N > \max(k, \ell)$ and note that the following statements hold:
\[\alpha^\epsilon(v_k) = \alpha_N^\epsilon(v_k),
\quad
\alpha^\epsilon(v_\ell) = \alpha_N^\epsilon(v_\ell)
.\]
Since $\alpha_N^\epsilon \in \A(G)$ it follows immediately that $\alpha^\epsilon$ is injective. If $v_\ell = \alpha_k^{-\epsilon}(v_k)$, then
$v_\ell = \alpha_N^{-\epsilon}(v_k)$ and so
$\alpha^\epsilon(v_\ell) = \alpha_N^\epsilon(v_\ell)
= \alpha_N^\epsilon(\alpha_N^{-\epsilon}(v_k))
= v_k$. Hence $\alpha^\epsilon$ is also surjective, and so we may be certain that $\alpha^\epsilon$ is a permutation of $V(G)$. 

For any (finite) ordered $n$-tuple $(w_1, \ldots, w_n)$ of elements in $V(G)$, there exists $N \in \N$ such that $\{w_1, \ldots, w_n\} \subseteq \{v_0, \ldots, v_N\}$, and so:
\begin{equation} \label{eq:finitary}
\alpha^\epsilon \left ( (w_1, \ldots, w_n) \right ) = \alpha_m^\epsilon \left ( (w_1, \ldots, w_n) \right ) \quad \quad \text{for all $m > N$.}
\end{equation}
It is now easy to see that $\alpha^\epsilon$ is an automorphism of $G$: the image of any edge or non-edge $\{w_1, w_2\}$ in $G$ under $\alpha^\epsilon$ is the same as its image under the automorphism
$\alpha_{N+1}^\epsilon$.

Finally, we wish to show that  $\alpha^\pi \neq \alpha^\epsilon$ for any pair $\epsilon, \pi$ of distinct elements in $E$.  For, let $k$ be the smallest index such that
$\pi_k\neq \epsilon_k$. We can assume that
$\varphi_k^{\epsilon_k} = \varphi_k$  and $\varphi_k^{\pi_k} = \id$.

We know that there is some $v \in F_{k+1}$ such that $v \neq \varphi_k(v) = \varphi_k^{\epsilon_k}(v)$. Now $v = v_\ell$ for some integer $\ell \geq 0$. Since $v_\ell \in F_\ell$, we must have $\ell > k$. Hence $\alpha^\epsilon(v) = \alpha_\ell^\epsilon(v) = \alpha_{k+1}^\epsilon(v)$.
Because  $v = \varphi_k^{\pi_k}(v)$ we infer
\begin{eqnarray*}\alpha^\epsilon(v) &=& (\varphi_0^{\epsilon_0}\varphi_1^{\epsilon_1} \cdots  \varphi_{k-1}^{\epsilon_{k-1}}) \varphi_{k}^{\epsilon_k}(v)\\
 &= &(\varphi_0^{\pi_0}\varphi_1^{\pi_1} \cdots \varphi_{k-1}^{\pi_{k-1}})\varphi_{k}^{\epsilon_k}(v)\\ & \neq & (\varphi_0^{\pi_0}\varphi_1^{\pi_1} \cdots \varphi_{k-1}^{\pi_{k-1}})\varphi^{\pi_k}(v)\\
 & =& \alpha^\pi(v)
 \end{eqnarray*}
This means that the $2^{\aleph_0}$ automorphisms $\alpha^\epsilon$, for $\epsilon \in E$, are all distinct.
\qed

We do not know how to extend Theorem \ref{thm:halin1} to graphs of higher cardinality although this may be feasible by the results of Kueker \cite{kueker} and Reyes \cite{Reyes} (see Section~\ref{sub:higher}).

\section{A bound on the distinguishing cost}

Given a graph with finite base it is not a priori clear whether it has finite 2-distinguishing cost. We show in Theorem \ref{thm:Cameron} that this is the case for graphs that have subdegree-finite automorphisms groups.

We need several lemmas for the  proof of Theorem \ref{thm:Cameron}. The first one is easy and well-known.

\begin{lemma}\label{le:finite vertex-stab}  If $G$ has a finite base and if $\A(G)$ is subdegree-finite, then all stabilizers of finite sets of vertices are finite.
\end{lemma}

\proof
Suppose $F$ is a finite set of vertices. Clearly $|\A(G)_F : \A(G)_{(F)}|$ is finite, and so $\A(G)_F$ is finite if and only if $\A(G)_{(F)}$ is finite.

Let $B$ be a finite base for $G$ and suppose $\A(G)$ is subdegree-finite.
For any vertex $v$, the set $C = \A(G)_v(B)$ is finite because all orbits of $\A(G)_v$ are finite.
Because $C$ contains $B$, it is a base for $G$. Hence, distinct automorphisms in $\A(G)_C$ cannot induce the same permutation of the finite set $C$, so $\A(G)_C$ must be finite.
Because $C$ is stabilized by $\A(G)_v$, it follows that $\A(G)_v$ is finite.
Therefore the point stabilizer $\A(G)_{(F)}$ is finite.
\qed

\begin{lemma}\label{le:infinite orbit}
Suppose $G$ is a graph with a finite base. If the automorphism group of $G$ is infinite, then it has an infinite orbit.
\end{lemma}

\proof Let $B$ be a finite base of $G$, and suppose all orbits of $\A(G)$ are finite.
Then $D:=\A(G)(B)$ is a finite set which is stabilized by $\A(G)$. Moreover, $D$ is a base for $G$ because it contains $B$.
Since no two distinct automorphisms in $\A(G)$ can induce the same permutation of $D$, it follows immediately that $\A(G)$ is finite.
\qed

The following two lemmas are similar to Lemma 2.1 and Theorem 3.4 (iii) of \cite{istw-2015}.

\begin{lemma}\label{subdegree} If $\A(G)$ is subdegree-finite and has an infinite orbit then, for any two finite sets of vertices  $Y, Z$ of $G$, there exists an element $\alpha\in \A(G)$ such that $Y\cap \alpha(Z)=\emptyset$.
 \end{lemma}

 \proof  Suppose a vertex $x$ lies in an infinite orbit of $\A(G)$.  Then there exists an infinite sequence $S=\{\alpha_i \in \A(G) : i \in\N \} \subseteq \A(G)$ such that
$\alpha_i^{-1} x = \alpha_j^{-1} x$
if and only if $i = j$.  Suppose there exist finite sets of vertices $Y$ and $Z$ such that $Y \cap \alpha_i (Z) \neq\emptyset$ for infinitely many $i\in\N$.  Then there exist $y \in Y$ and $z \in Z$ and an infinite subsequence $\{\alpha_{i_j} : j \in \N\} \subseteq S$ such that $\alpha_{i_j}z = y$. For each $j\in\N$ define $\beta_j := \alpha_{i_1}a_{i_j}^{-1}$, and notice that $\beta_j$ belongs to the stabilizer $\A(G)_y$. However, $\beta_j x = \beta_k x$  if and only if $j=k$, and so $\A(G)_y$ has an infinite orbit, a contradiction.
\qed

\begin{lemma}\label{le:infinite motion}  If $G$ has a finite base and $\A(G)$ is infinite and subdegree-finite, then $G$ has  infinite motion.
\end{lemma}

\proof
Let $G$ be a graph with a finite base $Y$ such that $\A(G)$ is infinite and subdegree-finite. Suppose $\beta \in \A(G)$ has finite motion. Then there exists a finite set $Z$ of vertices such that $\beta$ fixes no element in $Z$ and fixes (pointwise) every vertex not in $Z$. By Lemma~\ref{subdegree}, there is an $\alpha\in \A(G)$ such that  $\alpha(Z)\cap Y=\emptyset$.  Hence $\alpha\beta\alpha^{-1}$ fixes $Y$ pointwise, and so $\beta$ is the identity.
\qed

\begin{lemma}\label{le:main}
Suppose $G$ is a graph with a finite base such that $\A(G)$ is infinite and subdegree-finite. If $Y$ is a finite set of at least two vertices whose stabilizer is nontrivial, then there exist infinitely many vertices $v$ such that the stabilizer of $Y \cup \{v\}$ is a proper subgroup of the stabilizer of $Y$.
\end{lemma}

\proof Let $Y$ be a finite set of at least two vertices such that $\A(G)_Y$ is nontrivial, and let $D := \{\alpha \in \A(G) : \alpha(Y) \cap Y \not = \emptyset\}$. Suppose $D$ is infinite. Then there exists $y \in Y$ and an infinite subset $\{\alpha_i\}_{i \in \N} \subseteq D$ such that $\alpha_i(y) = \alpha_j(y)$ for all $i, j \in \N$. The elements in $\{\alpha_i^{-1}\alpha_1\}_{i \in \N}$ are pairwise distinct automorphisms which stabilize $y$, and so the stabilizer of $y$ is infinite. However, all stabilizers of finite sets are finite by Lemma~\ref{le:finite vertex-stab}. Hence $D$ must be finite.

Let $X := \bigcup_{\alpha \in D} \alpha(Y)$, and note that $X$ is a finite set which contains $Y$. By Lemma \ref{le:infinite motion}, $\A(G)$ has infinite motion, and therefore there exists a vertex $v \not\in X$ such that $\alpha(v) \neq v$ for some nontrivial $\alpha \in \A(G)_Y$. Note that we have infinitely many choices for $v$.

We show now that the stabilizer of $Y\cup \{v\}$ stabilizes $Y$. To see this, consider an automorphism $\gamma$ which stabilizes $Y\cup \{v\}$ but does not stabilize $Y$. Then $v \in \gamma(Y)$. The set $\gamma(Y) \cap Y$ is nonempty because $Y$ has at least two vertices, and thus $\gamma(Y) \subseteq X$. But then $v \in X$, contrary to its choice.

From this we infer that the stabilizer of $Y \cup \{v\}$ is contained in the stabilizer of $Y$, and they cannot be equal since $v$ is moved by the stabilizer of $Y$. \qed

\begin{theorem}\label{thm:Cameron} Suppose $G$ is a graph whose automorphism group is infinite and subdegree-finite. If $G$ has a finite base $B$, then $G$ has a finite distinguishing set with at most $\left\lceil\frac{5n}{2}\right\rceil - b(n) - 1$ elements, where $n = |B|$ and $b(n)$ denotes the number of 1s in the base-2 representation of $n$.
\end{theorem}

\proof Let $B$ be a finite base of $G$ of size $n$. If $n = 1$, then $G$ has a finite distinguishing set of size $1$ and the theorem holds. So, assume that $B$ contains at least two vertices.

Let $Y_0 := B$ and note that either $\A(G)_{Y_0}$ is trivial or, by Lemma~\ref{le:main}, we can find a vertex $v_1 \not \in Y_0$ such that the set $Y_1 := Y_0 \cup \{v_1\}$ satisfies $\A(G)_{Y_1} < \A(G)_{Y_0}$. Since $Y_1$ contains at least two vertices, we may repeat this process. Since $\A(G)_B$ is finite, we will eventually obtain a (possibly empty) finite set of vertices $\{v_1, \ldots, v_k\}$ such that the stabilizer of $Y_k = B \cup \{v_1, \ldots, v_k\}$ is trivial and
\[\{\id\} = \A(G)_{Y_k} < \cdots < \A(G)_{Y_1} < \A(G)_{Y_0} = \A(G)_{B}.  \]

Now we observe that $\A(G)_B$ is a  subgroup of the symmetric  group $\Sym(n)$ on $n$ elements.
By a theorem of Cameron, Solomon and Turull (\cite{casot-89}) the length of any chain of subgroups of $\Sym(n)$
 is bounded by  $\left\lceil\frac{3n}{2}\right\rceil - b(n) - 1$, where $b(n)$ denotes the number of 1s in the base-2 representation of $n$. Hence $k \leq \left\lceil\frac{3n}{2}\right\rceil - b(n) - 1$ and therefore $|Y_k| \leq \left\lceil\frac{5n}{2}\right\rceil - b(n) - 1$.
\qed

Despite the fact that we hold  subdegree-finiteness to be
 essential for Theorem~\ref{thm:Cameron}, we have no example of an infinite, connected graph with finite base and infinite distinguishing cost. We thus conclude this section with a question.

\begin{question} \label{question:highercardinals} Does there exist an infinite, connected graph with finite base and infinite distinguishing cost?
\end{question}

\section{Non-elementary results and methods}
\label{section:nonelementaryresults}

Halin's proof of Theorem~\ref{thm:halin} and our elementary proof of Theorem~\ref{thm:halin1} rely on nested sets of subgraphs, and these arguments are essentially topological. Such nested sets are used in the definition of the permutation topology, which is utilzed by Evans in \cite{Evans87}. We briefly describe this topology below.

Let $X$ be a countably infinite set and $\Sym(X)$ denote the symmetric group on $X$. 
Suppose $(X_i)_{i \geq 0}$ is a nested sequence of finite nonempty subsets of $X$ that covers $X$. In other words, $X_i \subset X_{i+1}$ and $ \bigcup_{i \geq 0}X_i = X$.

For distinct permutations $\alpha, \beta \in \Sym(X)$ one  defines the confluent of  $\alpha, \beta$ as
\[\conf(\alpha, \beta) = \min\{i \in \mathbb{N} \, | \, \exists x \in X_i: \alpha x \neq \beta x\}.\]
Hence,  the confluent is the maximum $i$ such that $\alpha$  and $\beta$  coincide on $X_i$,  and it is zero
if they 
differ on $X_0$.
Observe that $\conf$ depends on the choice of the sequence $X_i$.

Then we define a distance $d(\alpha, \beta)$ between $\alpha$ and $\beta$ by setting $d (\alpha, \beta) = 0$ for $\alpha = \beta$ and
$d (\alpha, \beta) = 2^{-\conf(\alpha, \beta)}$ otherwise. This is a well defined metric. In fact, it is an ultrametric (a metric where the triangle inequality has the form $d(\alpha, \gamma) \leq \max\{d(\alpha,\beta), d(\beta, \gamma)\}$).

All metrics of this kind define the same topology, the so-called \emph{permutation topology} on $\Sym(X)$,
under which
$\Sym(X)$ is a topological group, and it makes sense to speak of closed subgroups.

A \emph{relation of arity $n$} on $X$ is a set of $n$-tuples of $X$. A {\em finitary relational structure on $X$} is a pair $(X, R)$, where $R$ is a set of relations of finite arity on $X$. The automorphism group of a relational structure $(X, R)$ consists of those permutations $\alpha$ of $X$ which satisfy: $\alpha(r) = r$ for all $r \in R$. Graphs and digraphs are both examples of finitary relational structures.

It is well-known, and easy to prove, that $A \leq \Sym(X)$ is closed if and only if $A$ is the automorphism group of a finitary relational structure $\mathcal{R}$ on $X$ (see \cite[Theorem 2.6]{cameron:oligomorphic}, for example). Hence, if in the proof of Theorem~\ref{thm:halin1} we replace the word ``edge'' with the word ``relation'', $G$ with $\mathcal{R}$, $V(G)$ with $X$, and $\A(G)$ with $A$, then we obtain an elementary proof of the following theorem.

\begin{theorem} \label{thm:halin1b}
Let $X$ be a countable set and let $A$ be a closed subgroup of $\Sym(X)$. Then $A$ is countable or has cardinality $2^{\aleph_0}$, with $A$  countable if and only if $A$ has a finite base.
\end{theorem}

This is how far we can go with elementary methods. A more general result of this kind is the following theorem of Evans.

\begin{theorem}\label{thm:evans} {\rm  \cite[Theorem 1.1]{Evans87} (Evans 1987)}
Let $X$ be a countably infinite set and $G, H$ closed subgroups of $\Sym(X)$ with $H \leq G$. Then, either $|G:H| = 2^{\aleph_0}$ or $H$ contains the pointwise stabilizer in $G$ of some finite subset of $X$.
\end{theorem}

Theorem \ref{thm:halin1b}
is the special case of Theorem~\ref{thm:evans} for $H = \{\id\}$. The proof 
of Theorem~\ref{thm:evans} is not elementary. It uses
the fact that the coset space $(G:H)$ inherits a metric from the metric $d$ on $\Sym(X)$, 
then Evans
proceeds to construct uncountably many convergent series of automorphisms in the coset space if there is no finite subset of $X$ whose point stabilizer is contained in $H$.  The construction involves binary trees, and closure is needed to assure that the limits are not in $H$.

Evans' paper (\cite{Evans87}) also contains a short variant of this proof due to Theo Grundh\"{o}fer, which proves the slightly weaker assertion that $H$ contains the pointwise stabilizer in $G$ of some finite subset of $X$ if $|G:H| \leq \aleph_0$. It begins with the observation that $d^\ast(\alpha, \beta) = d(\alpha, \beta) + d(\alpha^{-1}, \beta^{-1})$ also is a metric on $\Sym(X)$, which is complete with respect to $d^\ast$. Moreover, open balls in both the $d$- and the $d^\ast$-topology are the same and translation by an element in $\Sym(X)$ is a homeomorphism of $(\Sym(X), d^\ast)$. If $|G:H| \leq \aleph_0$, then the right cosets of $H$ in $G$ decompose into a countable number of closed subsets. Now one observes that by Baire's Theorem the assumption that  a non-empty complete metric space is the countable union of closed sets implies that one of these closed sets has non-empty interior. Hence, some coset of $H$ has non-empty interior, and thus $H$ contains some open ball in $G$ around the identity. This, in turn, implies that $H$ contains the pointwise stabilizer in $G$ of some finite subset of $X$.

\subsection{Uncountable graphs}\label{sub:higher}
It would be tempting to extend Theorem~\ref{thm:halin1} to graphs with higher cardinalities. As we already mentioned  this might be feasible by the results of Kueker \cite{kueker} and Reyes \cite{Reyes}.

The only result we know of in this direction is by Dixon, Neumann and Thomas \cite{dix-86}. Not surprisingly, it needs  the Generalized Continuum Hypothesis (GCH).

\begin{theorem}\label{thm:dix0} {\rm   \cite[Theorem 2]{dix-86} (Dixon, Neumann and Thomas 1986)}
Let $X$ be an  infinite set of cardinality $\textswab{n}$ and $G$ a subgroup of $\Sym(X)$ with $|\Sym(X):G| < 2^{\textswab{n}}$. Then, under the assumption of the \emph{GCH}, there is a subset $Y$ of $X$ such that
$|Y| < \textswab{n}$ and $\Sym(X)_{(Y)} \leq G$.
\end{theorem}

This theorem is a corollary to their main theorem, which pertains to countable graphs, and does not need the Continuum Hypothesis.

\begin{theorem}\label{thm:dix} {\rm   \cite[Theorem 1]{dix-86} (Dixon, Neumann and Thomas 1986)}
Let $X$ be a countably infinite set and $G$ a subgroup of $\Sym(X)$ with $|\Sym(X):G| < 2^{\aleph_0}$.  Then there is a finite subset $F$ of $X$ such that
$$\Sym(X)_{(F)} \leq G \leq \Sym(X)_{\{F\}},$$ where $\Sym(X)_{(F)}$ is the point-stabilizer of $F$ and
$\Sym(X)_{\{F\}}$ the set stabilizer.
\end{theorem}

The proof of Theorem~\ref{thm:dix} uses \emph{moieties}, where a moiety is a subset $Y$ of $X$ such that $|Y| = |X\setminus Y| = |X|$, and a theorem of Sierpi\'nski which asserts that a countable set contains an almost disjoint family of $2^{\aleph_0}$ moieties, that is, the intersection of any two members of the family is finite.

Although the proof of  Theorem~\ref{thm:dix} extends to the setting of  Theorem~\ref{thm:dix0},  it does not seem to lend itself directly to a generalization of Theorem~\ref{thm:halin1}.

\end{document}